\documentclass[11pt]{article}
 \usepackage{amsmath,amssymb,amscd}
\newcommand{\pf}[1]{\trivlist \item[\hskip\labelsep\it #1\ ]}
\newcommand{\varpf}[1]{\trivlist \item[\hskip\labelsep\sc #1:]}
\newcommand{\qedbox}{$\rlap{$\sqcap$}\sqcup$}
\newcommand{\qed}{\qquad \qedbox \endtrivlist}
\newcommand{\varqed}{\hfill \rule{0.6em}{0.6em} \endtrivlist}
\newenvironment{proof}{\pf{Proof}}{\qed}

\newenvironment{examples}{\pf{Examples} 
   \begin{enumerate}}{\end{enumerate} \endtrivlist}
\newenvironment{items}{
  \begin{enumerate} 
                    
  }{\end{enumerate}}

\renewcommand{\iff}{\quad\Longleftrightarrow\quad}

\newcommand{\comps}{{\mathbb C}}
\newcommand{\reals}{{\mathbb R}}

\newcommand{\cstar}{{C^\ast}}

\newcommand{\hull}{{\operatorname{hull}}}
\newcommand{\Prim}{{\operatorname{Prim}}}

\newcommand{\supp}{{\operatorname{supp}}}

\newcommand{\cl}[1]{#1^-}

\newtheorem{theorem}{Theorem}
\newtheorem{lemma}[theorem]{Lemma}
\newtheorem{corollary}[theorem]{Corollary}
\newtheorem{proposition}[theorem]{Proposition}
\newtheorem{df}[theorem]{Definition}
\newenvironment{definition}{\begin{df} \rm}{\end{df}}  
\newtheorem{quest}{Question} 
\newenvironment{question}{\begin{quest} \rm}{\end{quest}}
\begin{document}
\title{An amenable, radical Banach algebra}
\author{\it Volker Runde}
\date{}
\maketitle
\begin{abstract}
We give an example of an amenable, radical Banach algebra, relying on results 
from non-abelian harmonic analysis due to H.\ Leptin, D.\ Poguntke and 
J.\ Boidol.
\end{abstract}
Let $A$ be a Banach algebra, and let $E$ be a Banach $A$-module. A bounded 
linear map $D \!: A \to E$ is called a {\it a derivation\/} if
\[
  D(ab) = a.Db + (Da).b \qquad (a,b  \in A).
\]
A derivation $D \!: A \to E$ is said to be {\it inner\/} if
\[
  Da = x.a -a.x \qquad (a \in A)
\]
for some $x \in E$. For any Banach $A$-module $E$, its dual space
$E^\ast$ is naturally equipped with a Banach $A$-module structure via
\[
  \langle x, a.\phi \rangle := \langle x.a ,\phi \rangle
  \quad\mbox{and}\quad 
  \langle x, \phi.a \rangle := \langle a.x, \phi \rangle
  \qquad (a \in A, \, x \in E, \, \phi \in E^\ast).
\]
\par
We can now give the definition of an amenable Banach algebra:
\begin{definition} \label{amB}
A Banach algebra $A$ is {\it amenable\/} if, for each Banach $A$-module $E$,
every derivation $D \!: A \to E^\ast$ is inner.
\end{definition}
\par
The notion of amenability for Banach algebras was introduced by B.\ E.\ 
Johnson in \cite{BEJ}. A locally compact group $G$ is called amenable if it 
possesses a translation-invariant mean, i.e.\ if there is a linear functional
$\phi \!: L^\infty(G) \to \comps$ satisfying
\[
  \phi(1) = \| \phi \| = 1 \qquad\mbox{and}\qquad 
  \phi(\delta_x \ast f) = \phi(f) \quad (x \in G, \, f \in L^\infty(G)).
\]
For instance, all abelian and all compact groups are amenable. For further 
information, see the monograph \cite{Pat}. In \cite{BEJ}, B.\ E.\ Johnson 
proved the following fundamental theorem which provides the motivation for
Definition \ref{amB}:
\begin{theorem}
Let $G$ be a locally compact group. Then $G$ is amenable if and only if 
$L^1(G)$ is amenable.
\end{theorem} 
\par
Since then, amenability has turned out to be an extremely fruitful concept in 
Banach algebra theory. We only would like to mention the following deep result
due to A.\ Connes (\cite{Connes}) and U.\ Haagerup (\cite{Haa}):
\begin{theorem}
A $\cstar$-algebra is amenable if and only if it is nuclear.
\end{theorem}
\par
In \cite{Rick}, P.\ C.\ Curtis asked the following question: 
\begin{question} \label{phil}
Is there an amenable, radical Banach algebra?
\end{question}
\par
In these notes, I would like to answer this question affirmativly. 
\par
Let $A$ be a Banach algebra, and let $\Prim(A)$ denote the space of its 
primitive ideals endowed with the Jacobson topology. Recall that a closed 
subset $F$ of $\Prim(A)$ is a {\it set of synthesis\/} for $A$ if $\ker(F)$ 
is the only closed ideal $I$ of $A$ such that $F = \hull(I)$. Otherwise, 
$F$ is called a set of non-synthesis for $A$
\begin{definition}
A Banach algebra is said to be {\it weakly Wiener\/} if the empty set is
a set of synthesis for $A$. 
\end{definition} 
\par
As is customary, we call a locally compact group $G$ weakly Wiener if
$L^1(G)$ is weakly Wiener. 
\par
The following proposition is easily verified:
\begin{proposition}
Suppose there is a locally compact group which is amenable, but not weakly
Wiener. Then there is an amenable, radical Banach algebra.
\end{proposition}
\begin{proof}
Since $G$ fails to be weakly Wiener, there is a proper, closed ideal $J$ of
$L^1(G)$ whose hull in $\Prim(L^1(G))$ is empty. Then $R := L^1(G)/J$ is a
radical Banach algebra which, being the quotient of an amenable Banach algebra,
has to be amenable.
\end{proof}
\par
It thus makes sense to ask:
\begin{question} \label{q2}
Is there a locally compact group which is amenable, but fails to be weakly 
Wiener?
\end{question}
\par
By the proposition, an affirmative answer to Question \ref{q2} entails an 
affirmative one to Question \ref{phil}. When I posed Question \ref{q2} to Jean 
Ludwig of Metz he was (to my surprise) not only able to answer the question 
immediately upon learning of it, but also claimed that Question \ref{q2} had 
already been settled in the eighties by J.\ Boidol. The example is the group
\[
  G_{4,9}(0) := \left\{ \left[\begin{array}{rrr}
                           e^t & x & e^t z \\ 0 & e^{-t} & y \\ 0 & 0 & 1
                        \end{array}\right] :
                        t,x,y,z \in \reals \right\}.
\]
Another way of describing $G_{4,9}(0)$ is as follows. The {\it Heisenberg 
group\/} is defined as
\[
  H_1 := \left\{ \left[\begin{array}{rrr}
                        1 & x & z \\ 0 & 1 & y \\ 0 & 0 & 1
                     \end{array}\right] :
                     x,y,z \in \reals \right\}.
\]
Let ${\mathfrak A}(H_1)$ denote the group of automorphisms of $H_1$, and define
$\phi \!: \reals \to {\mathfrak A}(H_1)$ through
\[
  \phi(t) \left( \left[\begin{array}{rrr}
                        1 & x & z \\ 0 & 1 & y \\ 0 & 0 & 1
                 \end{array}\right] \right) :=
          \left( \left[\begin{array}{rrr}
                  1 & e^{-t} x & z \\ 0 & 1 & e^t y \\ 0 & 0 & 1
          \end{array}\right] \right).
\]
Then we may identify $G_{4,9}(0)$ with the semidirect product $H_1 \times_\phi
\reals$. In particular, $G_{4,9}(0)$ is an abelian extension of the
nilpotent group $H_1$, thus solvable, and therefore amenable.
\par
What remains to be shown is that $G_{4,9}(0)$ is not weakly Wiener. As 
Ludwig claims, this was proved by Boidol. Apparently, Boidol never published
his finding. However, the proof of \cite[Theorem 6]{LP} can be modified
with the help of \cite[Lemma 1]{Boi} to yield the desired result.
\par
We require two lemmas, the first of which is completely elementary:
\begin{lemma} \label{l1}
Let $A$ be a Banach algebra which is weakly Wiener. Then every quotient of
$A$ is weakly Wiener.
\end{lemma}
\par
The second lemma is a variant of the lemma on \cite[p.\ 130]{LP}. For
a Banach algebra $A$ we write $M(A)$ to denote its double centralizer algebra.
\begin{lemma} \label{l2}
Let $A$ be a Banach algebra which is weakly Wiener, and let $p \in M(A)$
be an idempotent such that $ApA$ is dense in $A$. Then $pAp$
is weakly Wiener.
\end{lemma}
\begin{proof}
Let $J \subsetneq pAp$ be a closed ideal, and let $I$ denote the closed ideal
of $A$ generated by $J$. Obviously, $I = \cl{I_0}$ with
\[
  I_0 := J + AJ + JA + AJA,
\]
and consequently, $pIp \subset J$. This means that $I \subsetneq A$. Since,
by assumption, $A$ is weakly Wiener, there is a primitive ideal $P$ of $A$
such that $I \subset P$. Let $Q := pAp \cap P$. Assume that $Q = pAp$.
Then $(ApA)^2 \subset P$, and since $P$ is primitive, $ApA \subset P$.
The density of $ApA$ in $A$ yields $A =P$, which is a contradiction.
Let $\pi$ be an irreducible representation of $A$ on some linear space, 
say $E$, such that $\ker \pi = P$. Since $A$ is an ideal in $M(A)$, $\pi$
extends canonically to an irreducible representation of $M(A)$ on $E$, which
we denote by $\pi$ as well. Since $p A p \not\subset P$, we have $\pi(p) \neq
0$. Let $x \in \pi(p)E \setminus \{ 0 \}$. Then
\[
  \pi(pAp)x = \pi(pA)x = \pi(p)\pi(A)x = \pi(p)E,
\]
and consequently $(\pi |_{pAp}, \pi(p)E)$ is an irreducible representation
of $pAp$ (compare the proof of \cite[Theorem 26.14]{BD}), and $Q$ is primitive.
\end{proof}
\par
Furthermore, we require the theory of generalized $L^1$-algebras as
given in \cite{Lep}. Let $G$ be a locally compact group, and let $\cal A$
be a Banach $^\ast$-algebra with isometric involution such that $G$ acts on
$\cal A$ as a group of isometric $^\ast$-automorphisms; for $x \in G$, we 
write ${\cal A} \ni a \mapsto a^x$ for the automorphism implemented by
$x$. The Banach space $L^1(G,{\cal A})$ becomes a Banach $^\ast$-algebra with 
isometric involution via
\begin{eqnarray*}
  (f \ast g)(x) := \int_G f(xy)^{y^{-1}} g(y^{-1}) \, dy 
  \quad\mbox{and}\quad f^\ast(x) := \Delta_G(x)^{-1} (f(x^{-1})^x)^\ast & & \\
  (f,g \in L^1(G,{\cal A}), \, x \in G) & &
\end{eqnarray*}
where $dx$ denotes left Haar measure, and $\Delta_G$ is the modular function
on $G$.
\par
Our main result will be that for a specific choice of $\cal A$ the algebra 
$L^1(G,{\cal A})$ fails to be weakly Wiener and then use this to conclude that
$G_{4,9}(0)$ is not weakly Wiener.
\par
We shall be concerned with the following situation:
\begin{itemize}
\item $G$ and $H$ are locally compact groups, and
\item $G$ acts continuously and automorphically on $H$, i.e.\ there is a 
continuous mapping $G \times H \to H$, $(g,x) \mapsto x^g$ such that
$xy)^g = x^g y^g$, $(x^g)^h = x^{gh}$, and $x^1 = x$.
\end{itemize}
For each $g \in G$, $dx^g$ is left Haar measure. Thus, there is a positive
real number $\Delta_{G,H}(g)$ such that $dx^g = \Delta_{G,H}(g) \, dx$.
\par
For $f \in {\cal C}_0(H)$ and $h \in H$, we define $f_h$ and $f^h$ via
\[
  f^h(x) := f(hx) \quad\mbox{and}\quad f_h(x) := f(xh) \qquad
  (x \in H).
\]
We shall also require a subalgebra $\cal Q$ of ${\cal C}_0(H)$ with the 
following properties:
\begin{items}
\item $\cal Q$ is a $^\ast$-subalgebra of ${\cal C}_0(H)$ equipped with a 
Banach algebra $| \cdot |$ such that $|q^\ast| = |q| \geq \| q \|_\infty$.
\item For $q \in \cal Q$ and $h$, we have $q^h \in \cal Q$ and 
$|q^h| = |q|$.
\item For each $q \in \cal Q$, the map $H \ni h \mapsto q^h$ is continuous.
\item The compactly supported functions in $\cal Q$ form a dense subalgebra
${\cal Q}_{00}$. 
\item For every neighborhood $U$ of $1_H$, there is $u \in \cal Q$ such that
\begin{itemize}
\item[(a)] $u \not\equiv 0$ and $\supp(u) \subset U$,
\item[(b)] $u_h \in \cal Q$ for all $h \in H$, and
\item[(c)] the map $H \ni h \mapsto u_h$ is continuous.
\end{itemize}
\end{items}
\begin{examples}
\item Obviously, ${\cal C}_0(H)$ satisfies all the requirements.
\item In case $H$ is abelian, we may choose ${\cal Q} = A(H) \cong 
L^1(\hat{H})$.
\end{examples}
\par
Let $\cal Q$ be as described above. For $u \in \cal Q$ and $g \in G$, let
\[
  (u \circ g)(x) := u(x^{g^{-1}}) \qquad (x \in H).
\]
We assume moreover:
\begin{items}
\item[(vi)] For each $g \in G$, the map ${\cal Q} \ni q \mapsto
q \circ g$ is an isometric automorphisms of $\cal Q$.
\item[(vii)] For each $q \in \cal Q$, the map $G \ni g \mapsto q \circ g$ is
continuous.
\end{items}
\par
These assumptions are certainly true for the two examples given above. 
\par
In this situation, we may speak of $L^1(H,{\cal Q})$. For $f \in 
L^1(H,{\cal Q})$ and $g \in G$, let
\[
  f^g(x) := \Delta_{G,H}(g)^{-1} f(x^{g^{-1}}) \circ g \qquad (x \in H).
\]
Thus, for each $g \in G$, the mapping $L^1(H,{\cal Q}) \ni f \mapsto f^g$ is
an isometric $^\ast$-isomorphism, and we may speak of $L^1(G,L^1(H,{\cal Q}))$.
\par
We are finally in a position to state the main theorem of these notes:
\begin{theorem} \label{LPB}
Let $G$, $H$ and $\cal Q$ be given as above, let ${\cal A} :=
L^1(H,{\cal Q})$, and suppose that $\Delta_{G,H} \not\equiv 1$. Then
$L^1(G,{\cal A})$ is not weakly Wiener.
\end{theorem}
\begin{proof}
For the sake of brevity, write ${\cal L} := L^1(G,{\cal A})$. Let
${\mathfrak H} := L^2(H)$. We begin by defining a faithful $^\ast$-representation
$\rho$ of $\cal A$ on $\mathfrak H$. For $f \in \cal A$ and
$\xi \in \mathfrak H$, let --- note that we can view $f$ as a function on
$H \times H$ ---
\begin{equation} \label{eq1}
  (\rho(f)\xi)(x) = \int_H f(xy,y^{-1}) \xi(y^{-1}) \, dy 
  \qquad (x \in H).
\end{equation}
For $u,v \in {\cal Q}_{00}$, define $u \circ v \in \cal A$ through
\[
  (u \circ v)(x) := \Delta_H(x)^{1/2} u^x \bar{v}.
\]
Then we have
\begin{equation} \label{eq2}
  (u \circ v)(x,y) = \Delta(x)^{1/2} u(xy) \overline{v(x)} \qquad
  (x,y \in H).
\end{equation}
Letting $u'(x) := \Delta_H(x)^{1/2} u(x)$, we obtain from (\ref{eq1}) and
(\ref{eq2}) that
\begin{eqnarray*}
  (\rho(u \circ v) \xi)(x) & := & \int_H
  \Delta_H(xy^{-1})^{1/2} u(x) \overline{v(y^{-1})} \xi(y^{-1}) \, dy \\
  & = & u'(x) \int_H \xi(y^{-1}) \overline{v(y^{-1})} 
  \Delta_H(y^{-1})^{1/2}\, dy \\
  & = & \langle \xi, v' \rangle u'(x).
\end{eqnarray*}
Thus $\rho(u \circ v)$ is a rank one operator, which in case $\| u' \|_2 = 1$
is a projection.
\par
Fix a real valued function $u \in {\cal Q}_{00}$ such that $u(1_H) > 0$ and
$\| u' \|_2 = 1$, and let $p := u \circ u$. Since $\rho$ is faithful and
$\rho(p)$ has rank one, it follows that $p$ is a projection in $\cal A$
such that $p{\cal A}p = \comps p$. Define a projection $p^\# \in M({\cal L})$
by letting
\[
  (p^\# f)(g) := p^g f(g) \quad\mbox{and}\quad (f p^\#)(g) := f(g)p
  \qquad (g \in G).
\]
We wish to apply Lemma \ref{l2}. As is shown in the proof of 
\cite[Theorem 4]{LP}, ${\cal A}p{\cal A}$ is dense in $\cal A$, and as 
pointed out on \cite[p.\ 131]{LP} this implies the density of ${\cal L} p^\#
{\cal L}$ in $\cal L$. Thus, we are finished once we have established that
${\cal L}_p := p^\#{\cal L} p^\#$ fails to be weakly Wiener.
\par
Let $\cal K$ denote the compact linear operators on $\mathfrak H$. Again from
the proof of \cite[Theorem 4]{LP}, we see that $\rho({\cal A}) \subset \cal K$.
Define a unitary representation of $G$ on $\mathfrak H$ by letting
\[
  (\pi(g)\xi)(x) := \Delta_{G,H}(g)^{1/2} \xi(x^g) 
  \qquad (g \in G, \, x \in H, \, \xi \in {\mathfrak H}). 
\]
Then for $g \in G$, $a \in \cal A$, and $\xi \in \mathfrak H$, we have
\begin{eqnarray*}
   (\pi(g)^\ast\rho(a)\pi(g)\xi)(x) & = & \Delta_{G,H}^{1/2}(g)
   (\rho(a)\pi(g)\xi (x^{g^{-1}}) \\
   & = & \Delta_{G,H}^{1/2}(g) \int_H a(x^{g^{-1}}y,y^{-1})(\pi(g)\xi)(y^{-1}) 
   \, dy \\
   & = & \int_H  a(x^{g^{-1}}y,y^{-1})\xi((y^g)^{-1}) \, dy \\ 
   & = & \int_H \Delta_{G,H}(g)^{-1} a((xy)^{g^{-1}},(y^{-1})^{g^{-1}})
   \xi(y^{-1}) \, dy \\
   & = & (\rho(a^g)\xi)(x) \qquad (x \in H),
\end{eqnarray*}
i.e.\ $\pi$ implements the action of $G$ on $\rho(A)$. In what follows,
we shall suppress the symbol $\rho$ and view $\cal A$ as a subalgebra of
$\cal K$.
\par
Let $G$ act on $\cal K$ in the trivial way, and consider the generalized group
algebra $L^1(G,{\cal K})$. For $f \in \cal L$, let
\[
  f^\natural(g) := \pi(g)f(g) \qquad (g \in G).
\]
It is easy to see that $f^\natural \in L^1(G,{\cal K})$, and that
$\sigma \!: {\cal L} \to L^1(G,{\cal K}), \, f \mapsto f^\natural$ is a
faithful $^\ast$-ho\-mo\-mor\-phism.
We wish to compute $\sigma({\cal L}_p)$. For $f \in \cal L$ and $g \in G$, 
we have $(p^\# f p^\#)(g) = p^gf(g)p$, i.e.\ we have
\[
  f \in {\cal L}_p \iff \mbox{$f(g) = p^gf(g)p$ for almost all $g \in G$}.
\]
Thus, for $f \in \cal L$ and for almost all $g \in G$, we have
\[
  f^\natural(g) = \pi(g)f(g) = \pi(g)p^gf(g)p = p \pi(g)f(g)p =
  \phi(g)p,
\]
for some $\phi(g) \in \comps$. Let $| \cdot |_\ast$ denote the 
$\cstar$-norm on ${\cal K}$. Then we have
\[
  |\phi(g)| = |f^\natural(g)|_\ast = |f(g)|_\ast \leq \| f (g) \|.
\]
Consequently, $\sigma({\cal L}_p) \subset L^1(G) \cong p L^1(G,{\cal K})p$. 
View ${\cal Q}_{00}$ as
a subspace of $\mathfrak H$. Then cleary $\pi(G){\cal Q}_{00} \subset 
{\cal Q}_{00}$. In particular, we have
\[
  (\pi(g)u)(x) = \Delta_{G,H}^{1/2}(g)u(x^g) \qquad (g \in G, \, x \in H).
\]
Define
\[
  w(g) := \pi(g^{-1})u \circ u \qquad (g \in G).
\]
Is is easily seen (compare \cite[p.\ 129]{LP}) that 
\[
  p^g {\cal A}p = \comps w(g) \quad\mbox{and}\quad
  | \omega(g) |_\ast = 1.
\] 
It follows hat each $f \in {\cal L}_p$ has the form $f(g) = \phi(g)w(g)$.
Letting
\[
  \omega(g) := \| w(g) \|  \qquad (g \in G),
\]
we conclude that $\sigma({\cal L}_p)$ is the Beurling algebra 
$L^1(G,\omega)$ (it is shown on \cite[p.\ 130]{LP} that $\omega$ is indeed
a weight). For $g \in G$, we have
\begin{eqnarray*}
   \omega(g) & = & \| \pi(g^{-1})u \circ u \| \\
   & = & \int_H | (\pi(g^{-1})u)^x u | \Delta_H(x)^{1/2} \, dx \\
   & \geq & \int_H \| (\pi(g^{-1})u)^x u \|_\infty \Delta_H(x)^{1/2} \, dx \\
   & = & \int_H (\sup_{y \in H} |\pi(g^{-1})u(xy)u(y)|)
   \Delta_H(x)^{1/2} \, dx \\
   & \geq & \int_H |\Delta_{G,H}(g)^{1/2} u(1_H)u(x^{-1})|
   \Delta_H(x)^{1/2} \, dx \\
   & = &  \Delta_{G,H}(g)^{1/2} \Bigg[ 
   \underbrace{u (1_H) \int_H \Delta_H(x)^{1/2}|u(x^{-1})| 
   \,dx}_{:=\Omega>0}
   \Bigg].
\end{eqnarray*}
Assume that ${\cal L}_p$ is weakly Wiener. Then by Lemma \ref{l1},
$L^1(G,\omega)$ is weakly Wiener. Since $\Delta_{G,H}^{1/2}(\cdot)$ is
a homomorphism from $G$ into the abelian group $\reals \setminus \{ 0 \}$,
there is no loss of generality if we assume that $G$ is also abelian.
By assumption, there is $g \in G$, such that $\Delta_{G,H}(g)^{1/2} > 1$.
We thus have 
\[
  \sum_{n=1}^\infty \frac{\log\omega(g^n)}{n^2} \geq
  \sum_{n=1}^\infty \frac{\log\Omega}{n^2} + 
  \sum_{n=1}^\infty \frac{\log \Delta_{G,H}(g^n)^{1/2}}{n^2}
  = \sum_{n=1}^\infty \frac{\log\Omega}{n^2} + \frac{1}{2} 
  \sum_{n=1}^\infty  \frac{\Delta_{G,H}(g)}{n} =
  \infty.
\]
Thus, $L^1(G,\omega)$ does not satisfy the Beurling-Domar condition
(\cite[p.\ 132]{Rei}) for a Beurling algebra to be weakly Wiener, and
we have reached a contradiction (compare \cite{Boi}).
\end{proof}
\begin{corollary}
$G_{4,9}(0)$ is not weakly Wiener.
\end{corollary}
\begin{proof}
We identify $L^1(G_{4,9}(0))$ and $L^1(\reals,L^1(H_1))$. For $f \in
L^1(H_1)$, define 
\[
  \dot{f}(x,s) := \int f(x,y,z) e^{-i(sy+z)} \, dy \, dz.
\]
The map $L^1(H_1) \ni f \mapsto \dot{f}$ is easily seen to be an 
epimorphism onto ${\cal A} := L^1(\reals,{\cal Q})$ with ${\cal Q} :=
A(\reals)$, which in turn induces an epimorphism from
$L^1(G_{4,9}(0))$ onto $L^1(\reals,{\cal A})$. So, if $L^1(G_{4,9}(0))$
is weakly Wiener, the same must be true for $L^1(\reals,{\cal A})$ by
Lemma \ref{l1} (here, the action of $\reals$ on $\reals$ is given by
$(t,x) \mapsto e^{-t}x$). Consequently, $\Delta_{\reals,\reals}(t) =
e^{-t} \not \equiv 1$. From Theorem \ref{LPB}, we obtain that $L^1(\reals,
{\cal A})$ cannot be wekaly Wiener.
\end{proof}
\par
In view of the proposition, we finally obtain the promised answer to
Question \ref{phil}:
\begin{corollary}
There is an amenable, radical Banach algebra.
\end{corollary}
\par
It is not clear at all (and in fact extremely unlikely) that the amenable,
radical Banach algebra whose existence we have just proved is commutative.
Hence, the following question remains open (\cite[Problem 13]{Hel}):
\begin{question} 
Is there a commutative, amenable, radical Banach algebra?
\end{question}
\vfill
\begin{tabbing}
{\it Address\/}: \= Fachbereich 9 Mathematik \\
                 \> Universit\"at des Saarlandes \\
                 \> Postfach 151150 \\
                 \> 66041 Saarbr\"ucken \\
                 \> Germany \\ \medskip
{\it E-mail\/}:  \> {\tt runde@math.uni-sb.de}
\end{tabbing}

\begin{thebibliography}{0000}
%
\bibitem[Boi]{Boi} {\sc J.\ Boidol}, On a regularity condition for group 
algebras of non abelian locally compact groups. In {\sc N.\ Petridis},
{\sc S.\ K.\ Pichorides}, and {\sc N.\ Varopoulos} (ed.s), 
{\it Harmonic Analysis. Iraklion 1978\/}. Springer Verlag (1980), pp.\ 16--21.
%
\bibitem[B-D]{BD} {\sc F.\ F.\ Bonsall} and {\sc J.\ Duncan}, {\it Complete
Normed Algebras\/}. Springer Verlag (1973).
%
\bibitem[Con]{Connes} {\sc A.\ Connes}, On the cohomology of operator algebras.
{\it J.\ Funct.\ Anal.\/}\ {\bf 28\/} (1978), 248--253.
%
\bibitem[Hel]{Hel} {\sc A.\ Ya.\ Helemskii}, 31 problems of the homology of
the algebras of analysis. In {\sc V.\ P.\ Havin} and {\sc N.\ K.\ Nikolskii}
(ed.s), {\it Linear and Complex Analysis Problem Book 3. Part I\/}. Springer 
Verlag (1994), pp.\ 54--78.
%
\bibitem[Haa]{Haa} {\sc U.\ Haagerup}, All nuclear $\cstar$-algebras are 
amenable. {\it Invent.\ math.\/}\ {\bf 74\/} (1983), 305--319.
%
\bibitem[Joh]{BEJ} {\sc B.\ E.\ Johnson}, Cohomology in Banach algebras.
{\it Mem.\ Amer.\ Math.\ Soc.\/}\ {\bf 127\/} (1972).
%
\bibitem[Lep]{Lep} {\sc H.\ Leptin}, Verallgemeinerte $L^1$-Algebren und
projektive Darstellungen lokal kompakter Gruppen. {\it Invent.\ math.\/}\
{\bf 3\/} (1967), 257--281; {\bf 4\/} (1967), 68--86.
%
\bibitem[L-P]{LP} {\sc H.\ Leptin} and {\sc D.\ Poguntke}, Symmetry and 
nonsymmetry for locally compact groups. {\it J.\ Funct.\ Anal.\/}\ {\bf 33\/}
(1979), 119--134.
%
\bibitem[Loy]{Rick} {\sc R.\ J.\ Loy} (ed.), {\it Conference on Automatic
Continuity and Banach Algebras\/}. Austr.\ Nat.\ Univ.\ (1989).
%
\bibitem[Pat]{Pat} {\sc A.\ L.\ T.\ Paterson}, {\it Amenability\/}. American
Mathematical Society (1988).
%
\bibitem[Rei]{Rei} {\sc H.\ Reiter}, {\it Classical Harmonic Analysis and
Locally Compact Groups\/}. Oxford University Press (1968).
%
\end{thebibliography}
\end{document}